\documentclass{article}
\usepackage{graphics}
\usepackage[T1]{fontenc}
\usepackage{latexsym}
\usepackage{mathrsfs}

\title{A novel numerical technique used in the solution of ordinary differential equations with a mixture of integer and fractional derivatives}
\author{Jacek S. Leszczynski, Tomasz Blaszczyk \\
 Czestochowa University of Technology, \\ 
 Institute of Mathematics and Computer Science, \\
 ul. Dabrowskiego 73, 42-200 Czestochowa, Poland \\
 e-mail: \{jaclesz,tomblaszczyk\}@gmail.com}
\date{2007}

\begin{document}
\maketitle

\begin{abstract}
Using both fractional derivatives, defined in the Riemann-Liouville and Caputo senses, and classical derivatives of the
 integer order we examine different numerical approaches to ordinary differential equations. Generally we formulate
 some algorithms where four discrete forms of the Caputo derivative and three different numerical techniques of solving
 ordinary differential equations are proposed. We then illustrate how to introduce classical initial conditions into
 equations where the Riemann-Liouville derivative is included.
\end{abstract}

\section{Introduction}
In the past, fractional calculus was applied only from a mathematical point of view. The fundamental work was
done in~\cite{Miller,Oldham,Podlubny,Samko}. At present fractional calculus is extremely popular due to a rapid expansion in the
field of practical applications. Such applications have been used in physics and mechanics~\cite{Hilfer,Zaslavsky},
finance~\cite{Gorenflo2,Gorenflo1}, hydrology~\cite{Caputo01,Schumer} and many other disciplines.

Ordinary differential equations including a mixture of integer and fractional derivatives are a natural extension of
integer-order differential equations and give a novel approach to mathematical modeling many processes in nature. The
solution of a equation strongly depends on form of the equation and is still considered by many authors. It should be noted
that an analytical approach is limited to the linear form of equations and includes special functions such as  Fox and
Wright functions\cite{Gorenflo03,Kilbas} or the Mittag-Leffler function~\cite{Samko}. This greatly limits practical implementations,
i.e. sometimes it is very difficult to illustrate the solution in one simple chart. On the other hand, a numerical
solution~\cite{Blank,Diethlem03,Ford,Galucio,Leszcz00} is an alternative approach to  analytical one. However, this
approach has many disadvantages, i.e. the introduction of the initial conditions included in the Riemann-Liouville
derivative~\cite{Podlubny01}, the unreasonable assumption that a method applied to a single term equation is proper for solving
a multi-term equation~\cite{Diethlem04,ElSaka} etc. Against this background Ford~\cite{Ford} noticed that there can be
a~considerable gap between methods that perform well in theory and those whose implementations are effective.

In this paper we try to propose a numerical approach which will be more convenient in practical applications. We will give
a~numerical procedure for how to introduce classical initial conditions into an equation where the Riemann-Liouville derivative is
included. Here we will focus on such types of equation as
\begin{equation} \label{eq01}
 f\left( {x,y\left( x \right),D^1 y(x), \ldots, D^p y(x), D^{\alpha _1 } y\left( x \right),\ldots,D^{\alpha_m }
 y\left( x \right)} \right) = 0
\end{equation}
where y(x) is the solution obtained for the class of continuous functions, \\ $D^1y(x), \ldots, D^p y(x)$ are derivatives of the
integer order,\linebreak[4] $D^{\alpha_1}y(x), \ldots, D^{\alpha_m} y(x)$ are derivatives of the fractional order and\linebreak[4] 
$\alpha _1 ,\ldots,\alpha _m  \in R$ are real orders of a fractional derivative. We assume that the fractional derivative
is defined as the left-side Caputo derivative~\cite{Caputo00}
\begin{equation} \label{eq02}
 {}_{x_0 }^C D_x^\alpha  y\left( x \right) = \frac{1}{{\Gamma \left( {n - \alpha } \right)}}\int\limits_{x_0 }^x
 {\frac{{y^{\left( n \right)} \left( \tau  \right)}}{{\left( {x - \tau } \right)^{\alpha  - n + 1} }}d\tau }~~~~~~for~x>x_0
\end{equation}
and the left-side Riemman-Liouville derivative~\cite{Samko}
\begin{equation} \label{eq03}
 {}_{x_0 }D_x^\alpha  y\left( x \right) = \frac{1}{{\Gamma \left( {n - \alpha } \right)}}\frac{{d^n }}{{dx^n }}\int
 \limits_{x_0 }^x {\frac{{y\left( \tau  \right)}}{{\left( {x - \tau } \right)^{\alpha  - n + 1} }}d\tau }~~~~~~for~x>x_0 
\end{equation}
In above formulae, the notation $n=[\alpha]+1$ where $[\cdotp]$ is an integer part of a~real number. Moreover, we introduce
a~definition of the left-side Riemann-Liouville fractional integral~\cite{Oldham} as
\begin{equation} \label{eq04}
 {}_{x_0 }I_x^\beta  y\left( x \right) = \frac{1}{{\Gamma \left( \beta  \right)}}\int\limits_{x_0 }^x {\frac{{y\left( \tau
 \right)}}{{\left( {x - \tau } \right)^{1 - \beta } }}}d\tau~~~~~~for~x>x_0 
\end{equation}  
which will be used in our further calculations. Note that $\beta$ ($\beta > 0$) is the real order of Eqn.~(\ref{eq04}). On
the base of theory~\cite{Oldham} we use an expression
\begin{equation} \label{eq40}
 {}_{x_0 }^C D_x^\alpha  y\left( x \right) = {}_{x_0 }I_x^{n - \alpha } \left( {D^n y\left( x \right)} \right)
\end{equation}
which shows a relationship between the Caputo derivative~(\ref{eq02}) and the Riemann-Liouville integral~(\ref{eq03}). With
regard to papers~\cite{Blank,Diethlem00,Ford} in which numerical methods are used in the solution of fractional differential
equations, the authors mostly use the Caputo derivative. However, there is a small number of papers~\cite{Galucio} where
the authors use the Riemann-Liouville derivative. This small number of papers were confronted by the problem of how to introduce classical
initial conditions in the Riemann-Liouville derivative in order to obtain a~solution for a~class of continuous functions.
In this paper, we propose a way to avoid this problem. To be more precise, in every equation where the
Riemann-Liouville derivative occurs we will change it for the Caputo one. Following this we will discretize only the Caputo
derivative, except for one case where the real number of the Riemann-Liouville derivative dominates in the equation. In this case
we propose ${}_{x_0}D^{\alpha}_x y(x)=D^n {}_{x_0}I^{n-\alpha}_x y(x)$ and then discretize the left-side
Riemann-Liouville integral ${}_{x_0}I^{n-\alpha}_x$.

\section{Statement of the problem and its solution}

With regard to Eqn.~(\ref{eq01}) we limit our considerations to the equation which has the following form
\begin{equation} \label{eq42}
 D^p y\left( x \right) + \lambda \left\{
 \begin{array}{l}
  {}_{x_0 }^C D_x^\alpha  y\left( x \right) \\~~\\ 
  {}_{x_0 }D_x^\alpha  y\left( x \right) \\ 
 \end{array} \right. = 0
\end{equation}
where $p$ denotes an integer number being the derivatives order, \linebreak $\alpha \in \langle 0, 1)$ is the order of the fractional
derivative and $\lambda$ is an arbitrary real number. This simple form of the equation allows us to show how our methods work
properly in comparison to analytical solutions. Note that Eqn.~(\ref{eq42}) is the homogeneous ordinary differential
equation with a mixture of derivatives. The function $y(x)$ being the solution of this equation, strongly belongs to the
class of continuous functions. On the basis of our previous results~\cite{Leszcz00} we rewrite Eqn.~(\ref{eq42}) in an
explicit form. Consequently we obtain the three following types of equation:
\begin{itemize}
 \item $p>n$ for $p=2$, $\alpha\in\ \langle 0,1)$, $n=1$
       \begin{equation} \label{eq05}
        D^2 y\left( x \right) + \lambda {}_{x_0 }^C D_x^\alpha  y\left( x
        \right) = 0
       \end{equation}

       \begin{equation} \label{eq06}
        D^2 y\left( x \right) + \lambda {}_{x_0 } D_x^\alpha  y\left( x \right)
        = 0
       \end{equation} 
       It can be seen in above equations that the integer order of classical derivative dominates over the fractional one.
 \item $p=n$ for $p=1$, $\alpha\in\ \langle 0,1)$, $n=1$
       \begin{equation} \label{eq07}
        D^1 y\left( x \right) + \lambda {}_{x_0 }^C D_x^\alpha  y\left( x
        \right) = 0
       \end{equation}
       \begin{equation} \label{eq08}
        D^1 y\left( x \right) + \lambda {}_{x_0 } D_x^\alpha  y\left( x \right)
        = 0
       \end{equation}
       In this case we have equal integer orders for the classical and fractional derivative.
 \item $p<n$ for $p=0$, $\alpha\in\ \langle 0,1)$, $n=1$
       \begin{equation} \label{eq09}
        {}_{x_0 }^C D_x^\alpha  y\left( x \right) + \lambda y\left( x \right) =
        0
       \end{equation}
       \begin{equation} \label{eq10}
        {}_{x_0 } D_x^\alpha  y\left( x \right) + \lambda y\left( x \right) = 0
       \end{equation}
       The last case shows fractional ordinary differential equations which are well known in the literature. It may observe
       that the fractional order of the equation dominates over the integer one.
\end{itemize}

\subsection{Analytical solutions}

To compare our direct numerical results we are obligated to solve the above system of equations in an analytical way. In this
solution we will use a~general idea which transforms the Riemann-Liouville derivative to the Caputo one~\cite{Podlubny}.
Thus we have
\begin{equation} \label{eq11}
 {}_{x_0 }D_x^\alpha  y\left( x \right) = \sum\limits_{i = 0}^{n - 1} {\frac{{\left( {x - x_0 }
 \right)^{i-\alpha}}}{{\Gamma \left( {i - \alpha  + 1} \right)}}D^{i}y \left( {x_0 } \right) + {}_{x_0 }^C D_x^\alpha
 y\left( x \right)} 
\end{equation}
On the base of~\cite{Miller} we also use the Laplace transform. Following that the transform of the derivative of the
integer order $m\in N$ is
\begin{equation} \label{eq12}
 \mathcal{L} \left[ {D^{m}y \left( x \right)} \right] = s^m F\left( s \right) - \sum\limits_{k = 0}^{m - 1}
 {s^k D^{m - k - 1}y \left(x_0 \right)} 
\end{equation} 
The Laplace transform of the Caputo derivative is
\begin{equation} \label{eq13}
 \mathcal{L} \left[ {{}_{x_0 }^C D_x^\alpha  y\left( x \right)} \right] = s^\alpha  F\left( s \right) -
 \sum\limits_{k = 0}^{n - 1} {s^{\alpha  - k - 1} D^{k}y \left( x_0 \right)} 
\end{equation}
Using (\ref{eq11}) and (\ref{eq13}) we calculated the Laplace transform from the Riemann-Liouville derivative in the
following form
\begin{equation} \label{eq51}
 \begin{array}{c}
  \mathcal{L} \left[ {{}_{x_0 }D_x^\alpha  y\left( x \right)} \right] =\mathcal{L} \left[ {\sum\limits_{i = 0}^{n - 1}
  {\frac{{\left( {x - x_0 } \right)^{i - \alpha } }}{{\Gamma \left( {i - \alpha  + 1} \right)}}D^{i}y \left( {x_0 } \right)
  + {}_{x_0 }^C D_x^\alpha  y\left( x \right)} } \right] =  \\~~\\ 
  = \sum\limits_{i = 0}^{n - 1} {\frac{{D^{i}y \left( {x_0 } \right)}}{{\Gamma \left( {i - \alpha  + 1} \right)}}
  \frac{{\Gamma \left( {i - \alpha  + 1} \right)}}{{s^{i - \alpha  + 1} }}}  + s^\alpha  F\left( s \right) -
  \sum\limits_{j = 0}^{n - 1} {s^{\alpha  - j - 1} D^{j}y \left( {x_0 } \right)}  =  \\~~\\ 
  = s^\alpha  F\left( s \right) \\ 
 \end{array}
\end{equation}
It should be noted that Eqn.~(\ref{eq51}) is limited by initial conditions which are omitted here. In previous
considerations we assumed  function $y(x)$ to be continuous. Therefore Eqn.~(\ref{eq51}) is contrary to the Laplace
transform found in the literature~\cite{Miller} where initial conditions of non-integer order occur. This arises from an
assumption that function $y(x)$ is non-continuous.

In our analytical solutions we use also two additional transforms as
\begin{equation} \label{eq14}
 \mathcal{L}\left[ y \left( x \right) \right] = F\left( s \right)
\end{equation} 
\begin{equation} \label{eq15}
 \mathcal{L}\left[ {x^\alpha  } \right] = \frac{{\Gamma \left( {\alpha  + 1} \right)}}{{s^{\alpha  + 1} }}
\end{equation}

Using the above transforms in the set of equations (\ref{eq05})-(\ref{eq10}) and retransforming the results we obtain
analytical solutions. Including initial conditions
\begin{equation} \label{eq16}
 y \left( x_0 \right) = y_0,~~
 D^{1}y \left( x_0 \right) = \mathop {y'_0} 
\end{equation}
the analytical solution to Eqn.~(\ref{eq05}) is
\begin{equation} \label{eq17}
 \begin{array}{c}
  y\left( x \right) = y_0  + \mathop {y'_0 } \left( x-x_0 \right) E_{2 - \alpha ,2} \left( { - \lambda \left( x-x_0
  \right)^{2 - \alpha } } \right) \\~~\\ 
  D^{1}y \left( x \right) = \mathop {y'_0 } E_{2 - \alpha ,2} \left( { - \lambda \left( x-x_0 \right)^{2 - \alpha } } \right) \\
 + \mathop {y'_0 } \left( {2 - \alpha } \right) \left( x - x_0\right) E_{2 - \alpha ,2}^{\left( 1 \right)} \left( { - \lambda \left(
  x-x_0 \right)^{2 - \alpha } } \right) \\ 
 \end{array}
\end{equation}
where $ E_{\alpha ,\beta } \left( { - \lambda x^\alpha  } \right)$ denotes the Mittag-Leffler function~\cite{Samko} which is
defined as
\begin{equation} \label{eq18}
 E_{\alpha ,\beta } \left( { - \lambda x^\alpha  } \right) = \sum\limits_{i = 0}^\infty  {\frac{{\left( { - \lambda }
 \right)^i x^{\alpha i} }}{{\Gamma \left( {\alpha i + \beta } \right)}}}   
\end{equation}
It should be noted that 
$E_{\alpha ,\beta }^{\left( 1 \right)} \left( { - \lambda x^\alpha  } \right)$ occurs in solution (\ref{eq17}), which is the first derivative of the
Mittag-Leffler function~(\ref{eq18}) and is defined as
\begin{equation} \label{eq19}
 E_{\alpha ,\beta }^{\left( 1 \right)} \left( { - \lambda x^\alpha  } \right) = \sum\limits_{i = 1}^\infty  {\frac{{\left(
 { - \lambda } \right)^i ix^{\alpha i-1} }}{{\Gamma \left( {\alpha i + \beta } \right)}}}
\end{equation} 

Solving Eqn. (\ref{eq06}), where initial conditions (\ref{eq19}) are included, we have
\begin{equation} \label{eq20}
\begin{array}{c}
 y\left( x \right) = y_0 E_{2 - \alpha ,1} \left( { - \lambda \left( x-x_0 \right)^{2 - \alpha } } \right) \\
+ \mathop {y'_0 } \left( x-x_0 \right) E_{2 - \alpha ,2} \left( { - \lambda \left( x-x_0 \right)^{2 - \alpha } } \right) \\~~\\
 D^{1}y \left( x \right) = \left( {2 - \alpha } \right)y_0 E_{2 - \alpha ,1}^{\left( 1 \right)} \left( { - \lambda \left( x-x_0 \right)^{2 - \alpha } } \right)\\
~+ \mathop {y'_0 }  E_{2 - \alpha ,2} \left( { - \lambda \left( x-x_0 \right)^{2 - \alpha } } \right) \\ 
+ \left( {2 - \alpha } \right)\mathop {y'_0 } \left( x - x_0\right) E_{2 - \alpha ,2}^{\left( 1 \right)} \left( { - \lambda \left( x-x_0 \right)^{2 -
 \alpha } } \right) \\ 
 \end{array}
\end{equation}

Eqn.~(\ref{eq07}) with the initial condition $y(x_0)=y_0$ has the following solution
\begin{equation} \label{eq22}
 y\left( x \right) = y_0
\end{equation}

Eqn. (\ref{eq08}) with the initial condition $y(x_0)=y_0$ has the following solution
\begin{equation} \label{eq23}
 y\left( x \right) = y_0 E_{1 - \alpha ,1} \left( { - \lambda \left( x-x_0 \right)^{1 - \alpha } } \right) 
\end{equation}

Eqn. (\ref{eq09}) with the initial condition $y(x_0)=y_0$ has the following solution
\begin{equation} \label{eq24}
 y\left( x \right) = y_0 E_{\alpha ,1} \left( { - \lambda \left( x-x_0 \right)^{\alpha } } \right) 
\end{equation}

However Eqn. (\ref{eq10}) with the initial condition $y(x_0)=y_0$ has a trivial solution in the following form
\begin{equation} \label{eq25}
 y\left( x \right) =0 
\end{equation}
The above solution, arises from our assumption that we only consider a class of
continuous functions $y(x)$. In our next considerations we neglect this solution. It can be seen in literature~\cite{Ford,Galucio} that
the authors solve a similar type of equation to Eqn.~(\ref{eq10}), where an exciting function $f(x)$ is added on the right side of
the equation.

\subsection{Discrete forms of the Caputo derivative}

\begin{figure}[!h]
 \resizebox*{0.20\textwidth}{!}{\includegraphics{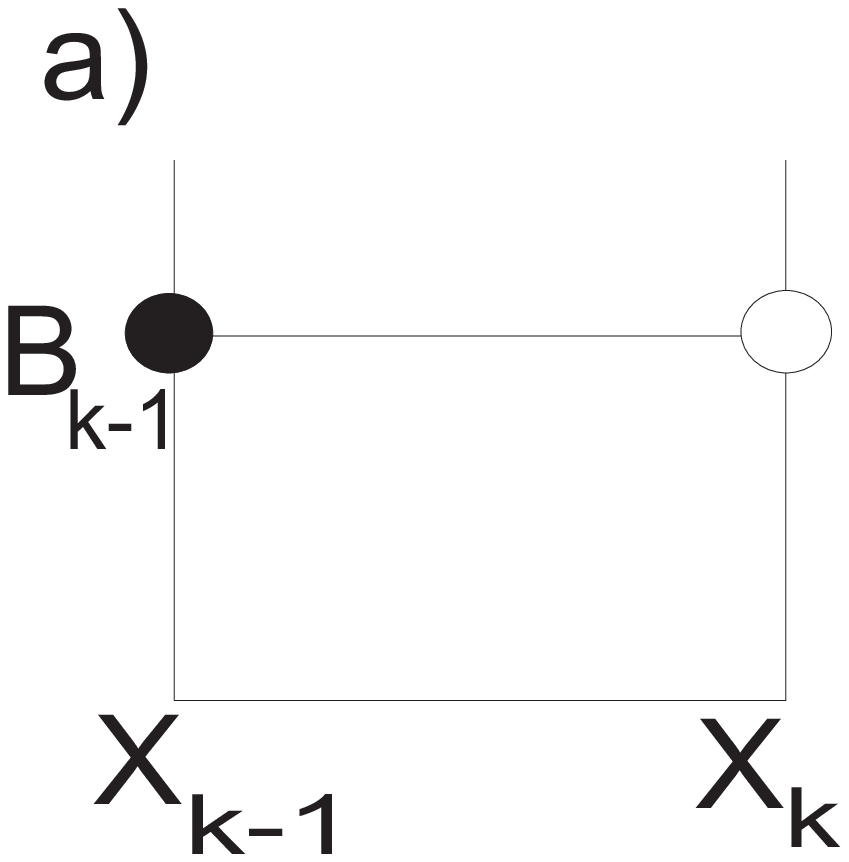}}
 \resizebox*{0.21\textwidth}{!}{\includegraphics{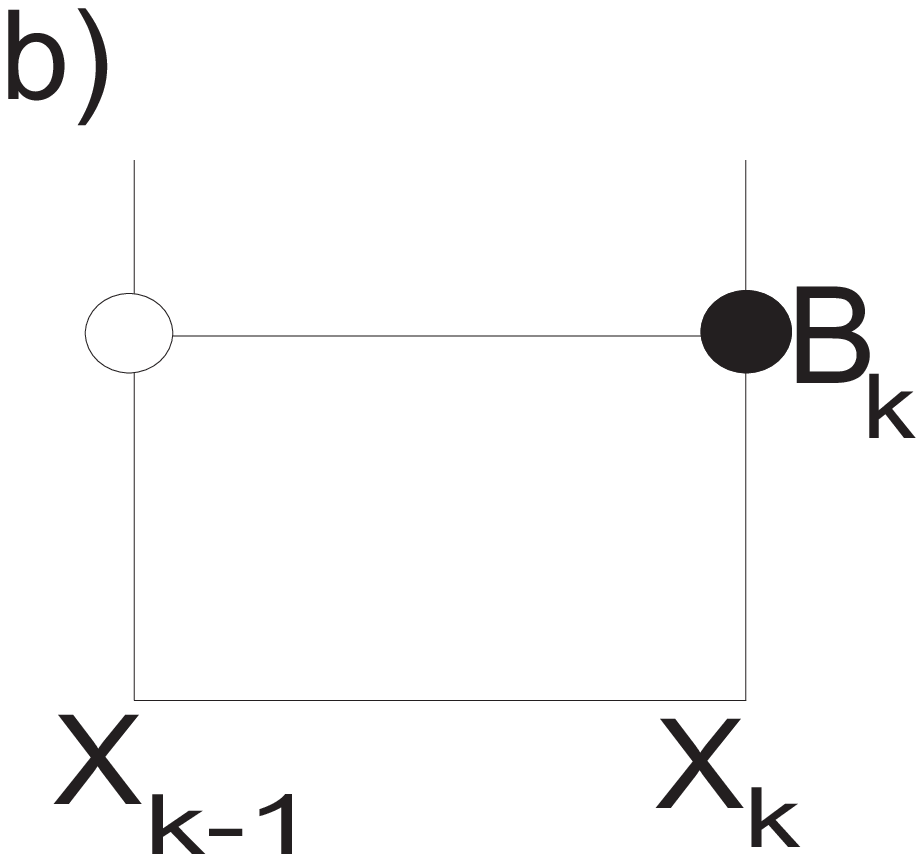}}
 \resizebox*{0.22\textwidth}{!}{\includegraphics{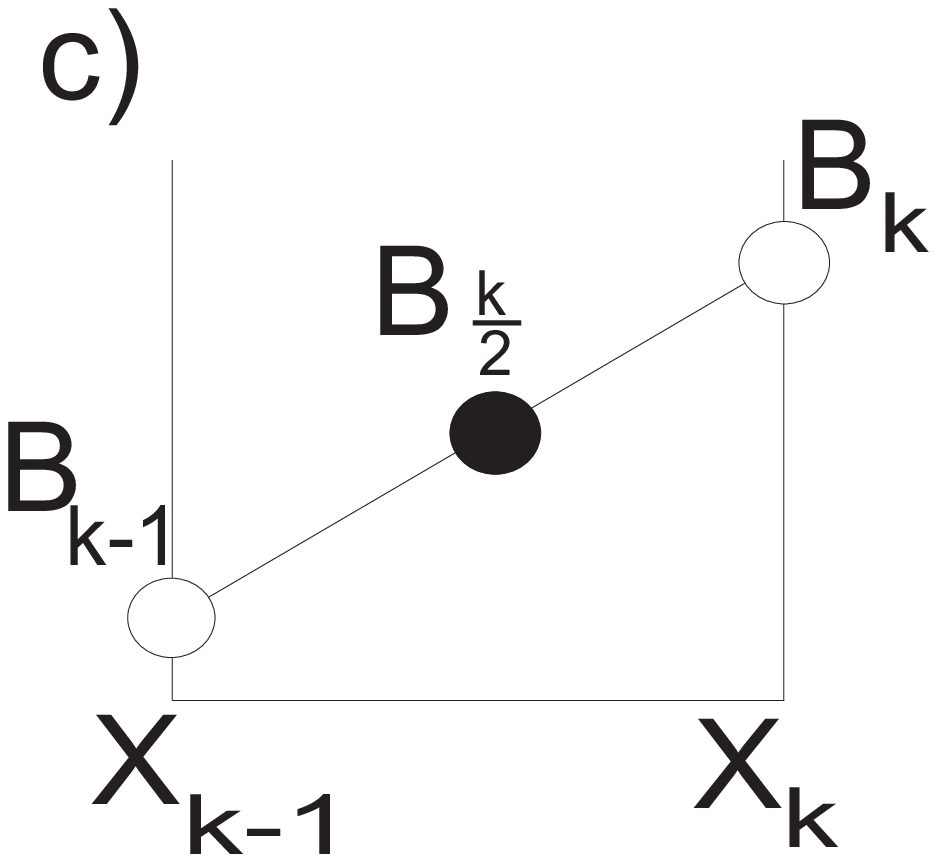}}
 \resizebox*{0.22\textwidth}{!}{\includegraphics{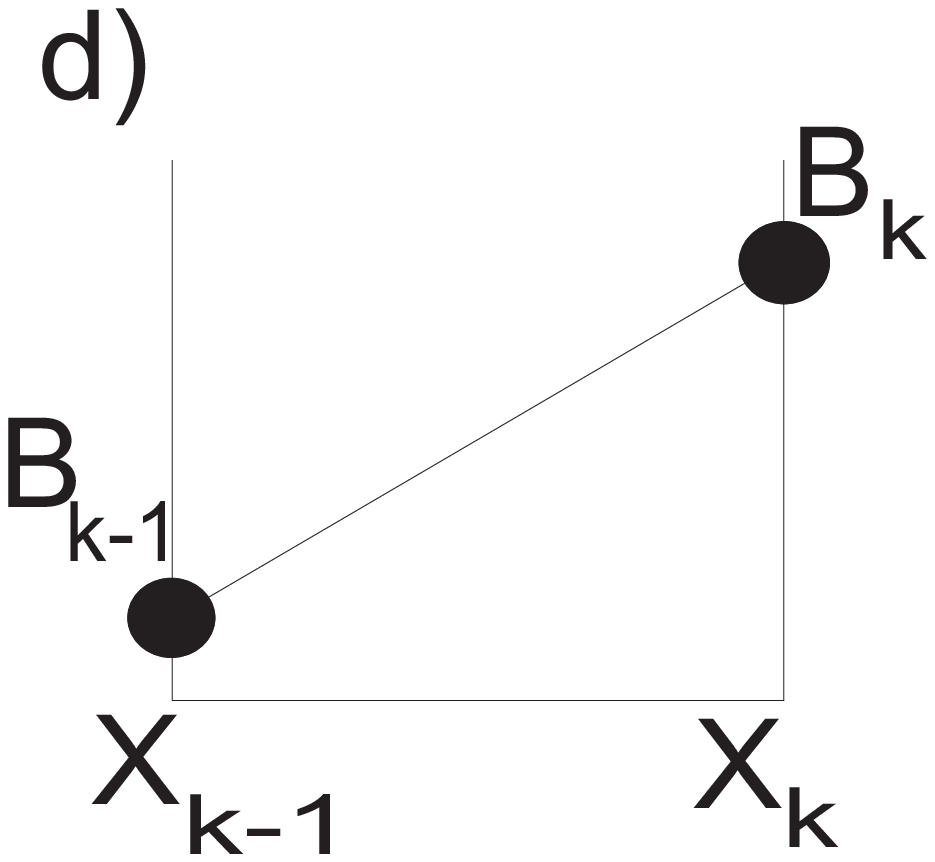}}
 \vspace{1cm}
\caption{}\label{fig01}
\end{figure}

There are many propositions on papers~\cite{Diethlem03,Gorenflo04,Oldham}, how to discretize fractional operators. Basically several discrete
forms are employed to take into account the different form of the function included in the fractional derivative. In this subsection we would
like to propose a general procedure for how to discretize the function. Let us consider an independent value $x$ which occurs
on a~length of calculations $\langle x_{0},x_{N}\rangle$, where $x_0$ and $x_N$ are the beginning and end of the range respectively. We
introduce a homogeneous grid $x_{0}<x_{1}<\ldots<x_{N}$. Fig.~\ref{fig01} \reversemarginpar \marginpar{\emph{Figure~\ref{fig01}}} shows four discrete forms of a~derivative
$D^n y(x)=B$ which is included in the Caputo derivative~(\ref{eq02}).

Taking into account the above discrete forms of the integer derivative we propose the following discrete schemes of the Caputo
derivative~(\ref{eq02}) as:
\begin{itemize}
\item the left-side form (case-I)
      \begin{equation} \label{eq26}
      \begin{array}{c}
      {}_{x_0 }^C D_x^\alpha  y\left( x \right) \cong \\
       \cong \frac{1}{{\Gamma \left(
       {n - \alpha  + 1} \right)}}\sum\limits_{k = 1}^N {B_{k - 1} \left[
       {\left( {x_N  - x_{k - 1} } \right)^{n - \alpha }  - \left( {x_N  - x_k }
       \right)^{n - \alpha } } \right]} 
       \end{array}      
       \end{equation}
      where $x \in \langle x_{k-1},x_{k}\rangle$ and
      $D^{n}y\left( x_{k-1} \right)=B_{k-1}$,
\item the right-side form (case-II)
      \begin{equation} \label{eq27}
      \begin{array}{c}
       {}_{x_0 }^C D_x^\alpha  y\left( x \right) \cong \\
       \cong \frac{1}{{\Gamma \left(
       {n - \alpha  + 1} \right)}}\sum\limits_{k = 1}^N {B_{k} \left[ {\left(
       {x_N  - x_{k - 1} } \right)^{n - \alpha }  - \left( {x_N  - x_k }
       \right)^{n - \alpha } } \right]} 
       \end{array}
       \end{equation}
      where $x \in \langle x_{k-1},x_{k}\rangle$ and
      $ D^{n}y\left( x_{k} \right)=B_{k}$,
\item the middle-side form (case-III)
      \begin{equation} \label{eq28}
      \begin{array}{c}
      {}_{x_0 }^C D_x^\alpha  y\left( x \right) \cong \\
      \cong \frac{1}{{\Gamma \left(
      {n - \alpha  + 1} \right)}}\sum\limits_{k = 1}^N {\frac{B_{k}+B_{k-1}}{2}
      \left[ {\left( {x_N  - x_{k - 1} } \right)^{n - \alpha }  - \left( {x_N
      - x_k } \right)^{n - \alpha } } \right]} 
      \end{array}
      \end{equation}
      where $x \in \langle x_{k-1},x_{k}\rangle$ and
      $\frac{D^{n}y\left( x_{k} \right)+D^{n}y\left( x_{k-1}
       \right)}{2}=\frac{B_{k}+B_{k-1}}{2}$. 
\item the linear form taken from~\cite{Leszcz01} (case-IV)
      \begin{equation} \label{eq29}
       \begin{array}{c}
        {}_{x_0 }^C D_x^\alpha  y\left( x \right) \cong \\
        \cong \frac{1}{{\Gamma \left(
        {n - \alpha } \right)}}\sum\limits_{k = 1}^N {\left\{ {\frac{{A_k }}{{n
        + 1 - \alpha }}\left[ {\left( {x_N  - x_k } \right)^{n + 1 - \alpha } –
        \left( {x_N  - x_{k - 1} } \right)^{n + 1 - \alpha } } \right]  }
        \right.}  \\ 
        + \left. {\frac{{A_k x_N  + B_k }}{{n - \alpha }}\left[
        {\left( {x_N  - x_{k - 1} } \right)^{n - \alpha }  - \left( {x_N  - x_k
        } \right)^{n - \alpha } } \right]} \right\} \\ 
       \end{array}
      \end{equation}
      where $A_{k}=\frac{D^{n}y\left( x_{k} \right)-D^{n}y\left( x_{k-1}
      \right)}{x_{k}-x_{k-1}}$, $B_{k}=D^{n}y\left( x_{k}\right)-A_{k}x_{k}$.
\end{itemize}

We try to use the above forms in numerical schemes to solve ordinary differential equations. It should be noted that
Eqn.~(\ref{eq26}) (case~I) is useful in the explicit scheme as, for example, the Euler's method~\cite{NumMet}. The next discrete
forms predicted by formulae (\ref{eq27}) (\ref{eq28}) and (\ref{eq29}) can use for any predictor-corrector
method~\cite{NumMet}.

Additionally, we also propose a discrete form of the \linebreak Riemann-Liouville fractional integral~(\ref{eq04}). We consider
a range \linebreak $\langle x_{k-1},x_{k}\rangle~\left( k=1,\ldots,N\right)$ for $k=1,\ldots,N$ and we assume a~constant value of function $y(x_k)=B_k$. In this case the integral operator~(\ref{eq04}) has the following form
\begin{equation} \label{eq30}
 {}_{x_0 }I_x^\beta  y\left( x \right) \cong \frac{1}{{\Gamma \left( {\beta  + 1} \right)}}\left[ {B_1 \left( 
 {x_N  - x_0 } \right) + \sum\limits_{k = 2}^N {\left( {B_k  - B_{k - 1} } \right)\left( {x_N  - x_{k - 1} }
 \right)^\beta  } } \right]
\end{equation}
where $B_{k}=y \left( x_{k} \right)$.

\subsection{Numerical methods of solving ordinary differential equations}

In this paper, we chose only three numerical methods which are used in the literature to solve an initial-value problem for
ordinary differential equations.

The Euler's method~\cite{NumMet} is an explicit one-step method. Using this method for any ordinary differential equation
of the first order we obtain
\begin{equation} \label{eq31}
 y_k  = y_{k - 1}  + hf\left( {x_{k - 1} ,y_{k - 1} } \right),~~k=1,\ldots,N
\end{equation}

The Adams method~\cite{NumMet} of the fourth order is the most popular method in the literature~\cite{Diethlem02,Diethlem05} which is used
to solve fractional differential equations. This is a~predictor-corrector method. It should be noted that for the method of
fourth order we have to determine three beginning values $y_1,y_2,y_3$ of the function. This can be done by using, for example,
the Euler's method. The Adams method for a~differential equation of the first order has the following form:
\begin{itemize}
 \item the predictor stage
       \begin{equation} \label{eq32}
        {}^{pr} y_k  = y_{k - 1}  + h\left( {\frac{{55}}{{24}}f_{k - 1}  -
        \frac{{59}}{{24}}f_{k - 2}  + \frac{{37}}{{24}}f_{k - 3}  -
        \frac{9}{{24}}f_{k - 4} } \right)
       \end{equation}
 \item the corrector stage
       \begin{equation} \label{eq33}
        {}^{cr} y_k  = y_{k - 1}  + h\left( {\frac{{19}}{{24}}f_{k - 1}  -
        \frac{5}{{24}}f_{k - 2}  + \frac{1}{{24}}f_{k - 3} } \right) +
        \frac{9}{{24}}hf_k 
       \end{equation}
       where $f_{k}  = f\left( {x_{k} ,{}^{pr}y_{k} } \right),~~f_{k - j}  =
       f\left( {x_{k - j} ,y_{k - j} } \right),~~j = 1,\ldots,4,\\
       k=4,\ldots,N$
\end{itemize}

The Gear's method~\cite{Gear} is also a predictor-corrector type method. This method has an advantage in cases where
only one call of the function is required in one step of the calculation. However, it needs higher derivatives in the algebraic form
at the beginning point $x_0$. This is a disadvantage of the method. The Gear's method for ordinary differential equations
of the first order is defined as
\begin{itemize}
 \item in the predictor stage
       \begin{equation} \label{eq34}
        {}^{pr}D^{i}y_k  = \sum\limits_{j = 0}^5 {P_{j - i} D^{j}y_{k - 1} } 
       \end{equation}
       where
       $P_l  = \left\{ \begin{array}{l} 0~~\mbox{for}~~l =  - 5,\ldots, -
       1 \\ \frac{{h^l }}{{l!}}~~\mbox{for}~~l = 0,\ldots,5 \\
       \end{array}\right.$
       and $i=0,\ldots,5$, $k=1,\ldots,N$
 \item in the corrector stage
       \begin{equation} \label{eq35}
        {}^{cr}D^{i}y_k  = {}^{pr} D^{i}y_k  + \frac{{i!}}{{2h^{i - 2} }}c_i
        \Delta {}^{cr}D^{2}y_k 
       \end{equation}
       where $\Delta {}^{cr}D^{2}y_k={}^{cr}D^{2}y_k-{}^{pr}D^{2}y_k$,
       ${}^{cr}D^{2}y_k=f\left( x_k, y_k, {}^{pr}D^{1}y_k \right) $ \\
       for $k=1,\ldots,N$, $i=0,\ldots,5$\\~~\\
       and coefficients are defined as\\
       $c_0 = \frac{3}{16}, c_1 = \frac{251}{360}, c_2 = 1, c_3 = \frac{11}{18},
       c_4 = \frac{1}{6}, c_5 = \frac{1}{60}$.\\~~\\
\end{itemize}
All the considered methods will be used in the next sections in order to construct algorithms.

\subsection{Algorithms}

In this subsection, we propose several algorithms which solve a~set of ordinary differential equations presented by general
formula~(\ref{eq42}).\\

\textbf{Algorithm 1.1}

At the beginning we consider Eqn.~(\ref{eq05}) with initial conditions (\ref{eq16}). On the base of the Euler's
method~(\ref{eq31}) and the left-side discrete form (\ref{eq26}) (case-I) of  Caputo derivative we then propose the
following algorithm
\begin{itemize}
 \item[\textbf{step 1}] Prediction of necessary data: initial conditions,
                        the fractional order $\alpha \in \langle 0, 1)$,
                        the total length of calculations 
                        $x\in \langle x_0, x_N \rangle$, the step of
                        calculations $h$.
 \item[\textbf{step 2}] Governing calculations: let $k = 1,\ldots, N$ then
               \begin{equation} \label{eq36}
                \begin{array}{l}
                 {}_{x_0 }^C D_x^\alpha  y_k = \frac{1}{{\Gamma \left( {2-\alpha } \right)}}\sum\limits_{j = 1}^k {D^1 y_{k - 1} \left[
                 {\left( {x_k  - x_{j - 1} } \right)^{1 - \alpha }  - \left(
                 {x_k - x_j } \right)^{1 - \alpha } } \right]}  \\~~\\ 
                  y_k = y_{k - 1} + hD^1 y_{k - 1} \\~~\\ 
                  D^1 y_k= D^1 y_{k-1} - h\lambda {}_{x_0 }^C D_x^\alpha  y_k \\ 
                \end{array}
               \end{equation}
\end{itemize}

\textbf{Algorithm 1.2}

Considering the same differential equation as presented in Algorithm~1.1 we present another approach. This is based on the
Adams method~(\ref{eq32}) (\ref{eq33}) and includes the linear-discrete form (\ref{eq29}) (case-IV) of the Caputo
derivative. Thus we have
\begin{itemize}
 \item[\textbf{step 1}] prediction of necessary data: initial conditions,
                        the fractional
                        order $\alpha \in \langle 0, 1)$, the total length of
                        calculations $x\in\langle x_0, x_N\rangle$, the step of
                        calculations h and additionally
                        ${}^{C}_{x_0}D^{\alpha}_x y_0 = 0$.
 \item[\textbf{step 2}] Introductory calculations: three initial values of the function
                        $y_1,y_2, y_3$ are calculated by the Euler's method
                        (algorithm~1.1).
 \item[\textbf{step 3}] Governing calculations: let $k = 4,\ldots, N$ then
               the predictor stage is
               \begin{equation} \label{eq43}
                \begin{array}{c}
                 {}^{pr}y_k  = y_{k - 1}  + h\left( {\frac{{55}}{{24}}D^1
                 y_{k-1}  - \frac{{59}}{{24}}D^1 y_{k - 2}  +
                 \frac{{37}}{{24}}D^1 y_{k - 3}  
                 - \frac{9}{{24}}D^1 y_{k - 4} } \right) \\~~\\ 
                 {}^{pr}D^1 y_k  = D^1 y_{k - 1}  \\
                 - h \lambda \left(
                 {\frac{{55}}{{24}}{}_{x_0 }^C D_x^\alpha  y_{k - 1}  -
                 \frac{{59}}{{24}}{}_{x_0 }^C D_x^\alpha  y_{k - 2} +
                 \frac{{37}}{{24}}{}_{x_0 }^C D_x^\alpha  y_{k - 3}
                -\frac{9}{{24}}{}_{x_0 }^C D_x^\alpha  y_{k - 4} } \right) \\ 
                \end{array}
               \end{equation}\\~~\\
               \begin{equation} \label{eq44}
                \begin{array}{c}
                 {}_{x_0 }^C D_x^\alpha  y_k  = \\
                 =\frac{1}{{\Gamma \left( { 
                 1-\alpha } \right)}}\sum\limits_{j = 1}^k \bigg\{ \frac{{D^1
                 y_j  - D^1 y_{j - 1} }}{{\left( {2 - \alpha } \right)h}}\left[
                 {\left( {x_k  - x_j } \right)^{2 - \alpha }  - \left( {x_k 
                 -x_{j - 1} } \right)^{2 - \alpha } } \right]     \\
                 +  \frac{{\left( {D^1 y_j  - D^1 y_{j - 1} }
                 \right)\left( {x_k  - x_j } \right) + hD^1 y_j }}{{\left( {
                 1-\alpha } \right)h}}\left[ {\left( {x_k  - x_{j - 1} }
                 \right)^{1 - \alpha }  - \left( {x_k  - x_j } \right)^{1 
                 -\alpha } } \right] \bigg\} \\ 
                \end{array}
               \end{equation}
               \\~~\\
               and the corrector stage is
               \begin{equation} \label{eq45}
                \begin{array}{c}
                 {}^{cr}y_k  = y_{k - 1}  + h\left( {\frac{{19}}{{24}}D^1 y_{k
                 -1}  - \frac{5}{{24}}D^1 y_{k - 2}  + \frac{1}{{24}}D^1 y_{k
                 -3} } \right) \\
                 + \frac{9}{{24}}hD^1 y_k  \\~~\\ 
                 {}^{cr}D^1 y_k  = D^1 y_{k - 1} \\
                 - h \lambda \left(
                 {\frac{{19}}{{24}}{}_{x_0 }^C D_x^\alpha  y_{k - 1}  -
                 \frac{5}{{24}}{}_{x_0 }^C D_x^\alpha  y_{k - 2}  +
                 \frac{1}{{24}}{}_{x_0 }^C D_x^\alpha  y_{k - 3} } \right)
                 -\frac{9}{{24}}h \lambda{}_{x_0 }^C D_x^\alpha  y_k  \\ 
                \end{array}
               \end{equation}
\end{itemize}
It should be noted that additional assumptions, as shown by step~2 and
step~3, need to be made for correct calculations to be achieved.\\

\textbf{Algorithm 1.3}

Still considering Eqn.~(\ref{eq05}) with initial conditions (\ref{eq16}) we can present another approach in comparison to
 previous ones. This uses the Gear's method (\ref{eq34}) (\ref{eq35}) and the middle-side discrete form (\ref{eq28})
(case-III) of  Caputo derivative.
\begin{itemize}
 \item[\textbf{step 1}] Prediction of necessary data: initial conditions,
                        the fractional order $\alpha\in\langle 0,1)$,
                        the total length of calculations
                        $x\in\langle x_0, x_N\rangle$,
                        the step of calculations $h$ and
                        additionally $D^2y_0=D^3y_0=D^4y_0=$ $=D^5y_0 = 0$.
 \item[\textbf{step 2}] Governing calculations: let $k=1,\ldots, N$ then
               Eqn.~(\ref{eq34}) is the predictor stage and beginning of the
               correction stage is
               \begin{equation} \label{eq49}
                \begin{array}{l}
                 {}_{x_0 }^C D_x^\alpha  y_k  = \frac{1}{{\Gamma \left( {2
                 -\alpha } \right)}}\sum\limits_{j = 1}^k {\frac{{D^1 y_k  + D^1
                 y_{k - 1} }}{2}\left[ {\left( {x_k  - x_{j - 1} } \right)^{1 
                 -\alpha }  - \left( {x_k  - x_j } \right)^{1 - \alpha } }
                 \right]}\\~~\\
                 \Delta {}^{cr}D^{2}y_k=-\lambda{}_{x_0 }^C D_x^\alpha  y_k  -
                 {}^{pr}D^{2}y_k
                \end{array} 
               \end{equation}
               and it follows that we apply the correction stage given by
               Eqn.~(\ref{eq35}).
\end{itemize}

\textbf{Algorithm 2.1}

Now let us consider Eqn. (\ref{eq06}) with initial conditions (\ref{eq16}). It should be noted that this type of equation includes the
Riemann-Liouville derivative. We apply the Euler's method~(\ref{eq31}) and Eqn. (\ref{eq11}) which transforms the
Riemann-Liouville derivative to the Caputo one. It should be remembered that Eqn.~(\ref{eq11}) includes a~set of initial
conditions in the general form
$\sum\limits_{i = 0}^{n - 1} {\frac{{\left( {x - x_0 } \right)^{i-\alpha}}}{{\Gamma \left( {i - \alpha  + 1}
 \right)}}}D^{i}y \left( {x_0 } \right)$.
In the case of the lower limit $x$ tends to $x_0$ and then the above expression is infinite. This means that we cannot use a full numerical approach in order to solve the class of ordinary differential equations where the
Riemann-Liouville derivative is included.  However an assumption of homogeneous initial conditions i.e.
$D^I y\left(x_0\right)=0$ avoids the problem. Eqn.~(\ref{eq06}) has a correct analytical solution
(\ref{eq20}) and the function $y(x)$ has a~finite value at the beginning point $x_0$ which is contrary to previous
considerations. With regard to our assumption that function $y(x)$ belongs to the class of continuous functions we are
obligated to improve the factor $\left(x-x_0\right)^{i-\alpha}$ included in Eqn. (\ref{eq11}) in order to restrict the
function continuity. Therefore we put $\left(x^{*}-x_0\right)^{i-\alpha}$ into Eqn. (\ref{eq11}) instead of
$\left(x-x_0\right)^{i-\alpha}$ where $x^{*}=x+c$ and $c = \left( \Gamma \left( 1-\alpha \right) \right)^{- \alpha^{-1}}$.
Point $c$ is where the shift from initial conditions to homogeneous initial conditions occurs. Our assumption allows for the correct behaviour of factor $\left(x^{*}-x_0\right)^{i-\alpha}=\left(x+c-x_0\right)^{i-\alpha}$ in Eqn. (\ref{eq11}) in
the class of continuous functions especially when $x$ tends to $x_0$. Taking into account the above considerations we have the
following algorithm
\begin{itemize}
 \item[\textbf{step 1}] prediction of necessary data: initial conditions, the fractional
               order \linebreak[4]$\alpha\in\langle 0,1)$, the total length of calculations
               $x\in\langle x_0, x_N\rangle$, the step of calculations $h$.
 \item[\textbf{step 2}] Governing calculations: let $k = 1,\ldots, N$ then
               \begin{equation}
                \begin{array}{c} \label{eq37}
                 {}_{x_0 }D_x^\alpha  y_k = \frac{{\left( {x_{k - 1}^*  - x_0 }
                 \right)^{ - \alpha } }}{{\Gamma \left( {1 - \alpha }
                 \right)}}y_0   \\
                 +\frac{1}{{\Gamma \left( {2 - \alpha }
                 \right)}}\sum\limits_{j = 1}^k {D^1 y_{k - 1} \left[ {\left(
                 {x_k  - x_{j - 1} } \right)^{1 - \alpha }  - \left( {x_k  - x_j
                 } \right)^{1 - \alpha } } \right]}  \\~~\\ 
                 y_k = y_{k - 1} + hD^1 y_{k - 1} \\~~\\ 
                 D^1 y_k = D^1 y_{k - 1} - h\lambda {}_{x_0 }D_x^\alpha  y_k \\ 
                \end{array}
               \end{equation} 
               where $x_{k - 1}^*  = x_{k - 1}  + hc$.
\end{itemize}

\textbf{Algorithm 2.2}

The Euler's method has been used in algorithm 2.1. Considering Eqn.~(\ref{eq06}) with initial conditions (\ref{eq16})
and using the Adams method (\ref{eq32}) with the linear-discrete form (\ref{eq29}) (case-IV) of  Caputo derivative we
obtain
\begin{itemize}
 \item[\textbf{step 1}] prediction of necessary data: initial conditions,
                        the fractional order $\alpha\in\langle 0,1)$,
                        the total length of calculations
                        $x\in\langle x_0, x_N\rangle$,
                        the step of calculations $h$,
                        ${}^{C}_{x_0}D^{\alpha}_x y_0 = 0$
                        and $c = \left( \Gamma \left(
                        1-\alpha \right) \right)^{- \alpha^{-1}}$.
 \item[\textbf{step 2}] Introductory calculations: three initial values of the function
               $y_1,y_2, y_3$ are calculated by the Euler's method
               (algorithm~2.1).
 \item[\textbf{step 3}] Governing calculations: let $k = 4,\ldots, N$ then
               the predictor stage is
               \begin{equation} \label{eq46}
                \begin{array}{c}
                 {}^{pr}y_k  = y_{k - 1}  + h\left( {\frac{{55}}{{24}}D^1 y_{k
                 -1}  - \frac{{59}}{{24}}D^1 y_{k - 2}  + \frac{{37}}{{24}}D^1
                 y_{k - 3}  - \frac{9}{{24}}D^1 y_{k - 4} } \right) \\~~\\  
                 {}^{pr}D^1 y_k  = D^1 y_{k - 1}  \\
                 - h\lambda\left(
                 {\frac{{55}}{{24}}{}_{x_0 } D_x^\alpha  y_{k - 1}  -
                 \frac{{59}}{{24}}{}_{x_0 } D_x^\alpha  y_{k - 2}  +
                 \frac{{37}}{{24}}{}_{x_0 } D_x^\alpha  y_{k - 3}  -
                 \frac{9}{{24}}{}_{x_0 } D_x^\alpha  y_{k - 4} } \right) \\ 
                \end{array}
               \end{equation}

               \begin{equation} \label{eq47}
                \begin{array}{c}
                 {}_{x_0 } D_x^\alpha  y_k  =\frac{{\left( {x_{k - 1}^*  - x_0 }
                 \right)^{ - \alpha } }}{{\Gamma \left( {1 - \alpha }
                 \right)}}y_0 \\
                 + \frac{1}{{\Gamma \left( {1 - \alpha }
                 \right)}}\sum\limits_{j = 1}^k \bigg\{ \frac{{D^1 y_j  - D^1
                 y_{j - 1} }}{{\left( {2 - \alpha } \right)h}}\left[ {\left(
                 {x_k  - x_j } \right)^{2 - \alpha }  - \left( {x_k  - x_{j - 1}
                 } \right)^{2 - \alpha } } \right]   \\ 
                 + \frac{{\left( {D^1 y_j  - D^1 y_{j - 1} }
                 \right)\left( {x_k  - x_j } \right) + hD^1 y_j }}{{\left( {1
                 -\alpha } \right)h}}\left[ {\left( {x_k  - x_{j - 1} }
                 \right)^{1 - \alpha }  - \left( {x_k  - x_j } \right)^{1 
                 -\alpha } } \right] \bigg\} \\ 
                \end{array}
               \end{equation}
               and the corrector stage is
               \begin{equation} \label{eq48}
                \begin{array}{c}
                 {}^{cr}y_k  = y_{k - 1}  + h\left( {\frac{{19}}{{24}}D^1 y_{k
                 -1}  - \frac{5}{{24}}D^1 y_{k - 2}  + \frac{1}{{24}}D^1 y_{k
                 -3} } \right) + \frac{9}{{24}}hD^1 y_k  \\~~\\ 
                 {}^{cr}D^1 y_k  = D^1 y_{k - 1}  +\\
                 - h\lambda\left(
                 {\frac{{19}}{{24}}{}_{x_0 } D_x^\alpha  y_{k - 1}  -
                 \frac{5}{{24}}{}_{x_0 } D_x^\alpha  y_{k - 2}  +
                 \frac{1}{{24}}{}_{x_0 } D_x^\alpha  y_{k - 3} } \right)
                 -\frac{9}{{24}}h\lambda{}_{x_0 } D_x^\alpha  y_k  \\ 
                \end{array}
               \end{equation}
\end{itemize}

\textbf{Algorithm 2.3}

We present another algorithm on the base of Algorithms 2.1 and 2.2. We still consider Eqn. (\ref{eq06}) with the initial conditions (\ref{eq16}) and then we use the Gear's method (\ref{eq34}) (\ref{eq35}) with the middle-side discrete form (\ref{eq28}) (case-III) of  Caputo derivative. We have
\begin{itemize}
 \item[\textbf{step 1}] prediction of necessary data: initial conditions,
                        the fractional order $\alpha\in\langle 0,1)$,
                        the total length of calculations
                        $x\in\langle x_0, x_N\rangle$,
                        the step of calculations $h$,
                        $c = \left( \Gamma \left(1-\alpha \right) 
                        \right)^{- \alpha^{-1}}$ and 
                        $D^2y_0=D^3y_0=$ $=D^4y_0=D^5y_0 = 0$.
 \item[\textbf{step 2}] Governing calculations: let $k=1,\ldots, N$ then
               Eqn.~(\ref{eq34}) is the predictor stage and beginning of the
               correction stage is
               \begin{equation} \label{eq50}
                \begin{array}{c}
                 {}_{x_0 } D_x^\alpha  y_k  = \frac{1}{{\Gamma \left( {2
                 -\alpha } \right)}}\sum\limits_{j = 1}^k {\frac{{D^1 y_k  + D^1
                 y_{k - 1} }}{2}\left[ {\left( {x_k  - x_{j - 1} } \right)^{1 
                 -\alpha }  - \left( {x_k  - x_j } \right)^{1 - \alpha } }
                 \right]}\\~~\\
                 \Delta {}^{cr}D^{2}y_k=-\lambda{}_{x_0 } D_x^\alpha  y_k  -
                 {}^{pr}D^{2}y_k
                \end{array} 
               \end{equation}\\
               and it follows that we apply the correction stage given by
               Eqn.~(\ref{eq35}).
\end{itemize}

\textbf{Algorithm 3.1}

Next we consider another class of ordinary differential equations given by formula (\ref{eq07}) with the initial condition $y(x_0)=y_0$. Using the Euler's method~(\ref{eq31}) and the left-side discrete form (\ref{eq26}) (case-I) of the Caputo derivative we then propose the following algorithm
\begin{itemize}
 \item[\textbf{step 1}] prediction of necessary data: initial conditions, the fractional
               order $\alpha \in \langle 0, 1)$, the total length of
               calculations $x\in \langle x_0, x_N \rangle$, the step of
               calculations $h$ and $D^{1} y_0 = 0$.
 \item[\textbf{step 2}] Governing calculations: let $k = 1,\ldots, N$ then
               \begin{equation} \label{eq38}
                \begin{array}{c}
                 {}_{x_0 }^C D_x^\alpha  y_k = \frac{1}{{\Gamma \left( {2
                 -\alpha } \right)}}\sum\limits_{j = 1}^k {D^1 y_{k - 1} \left[
                 {\left( {x_k  - x_{j - 1} } \right)^{1 - \alpha }  - \left(
                 {x_k  - x_j } \right)^{1 - \alpha } } \right]}  \\~~\\ 
                 y_k = y_{k - 1} - h \lambda {}_{x_0 }^C D_x^{\alpha} y_{k - 1}
                 \\~~\\ D^1 y_k = - \lambda {}_{x_0 }^C D_x^\alpha  y_k \\
                \end{array}
               \end{equation}
\end{itemize}

To simplify our considerations we neglect other algorithms which can be applied to solve Eqn. (\ref{eq07}). On the basis of
previous explanations of predictor-corrector methods one can construct interesting procedures self-reliantly.\\

\textbf{Algorithm 4.1}

Using the Euler's method (\ref{eq31}) and the left-side discrete form (\ref{eq26}) (case-I) of the Caputo derivative
we construct an algorithm which solves Eqn. (\ref{eq08}) with the initial condition $y\left(x_0\right)=y_0$. Thus we
obtain
\begin{itemize}
 \item[\textbf{step 1}] prediction of necessary data: initial conditions, the fractional
               order $\alpha \in \langle 0, 1)$, the total length of
               calculations $x\in \langle x_0, x_N \rangle$, the step of
               calculations $h$, $D^{1} y_0 = 0$ and $c = \left( \Gamma \left(1-
               \alpha \right) \right)^{- \alpha^{-1}}$.
 \item[\textbf{step 2}] Governing calculations: let $k = 1,\ldots, N$ then
               \begin{equation} \label{eq39}
                \begin{array}{c}
                 {}_{x_0 } D_x^\alpha  y_k =\frac{{\left( {x_{k - 1}^*  - x_0 }
                 \right)^{ - \alpha } }}{{\Gamma \left( {1 - \alpha }
                 \right)}}y_0   \\
                 + \frac{1}{{\Gamma \left( {2 - \alpha }
                 \right)}}\sum\limits_{j = 1}^k {D^1 y_{k - 1} \left[ {\left(
                 {x_k  - x_{j - 1} } \right)^{1 - \alpha }  - \left( {x_k-x_j}
                 \right)^{1 - \alpha } } \right]}  \\~~\\ 
                 y_k = y_{k - 1}-h \lambda {}_{x_0} D_x^{\alpha} y_{k-1} \\~~\\
                 D^1 y_k = - \lambda {}_{x_0 } D_x^\alpha  y_k \\
                \end{array}
               \end{equation}
               where $x_{k - 1}^*  = x_{k - 1}  + hc$.
\end{itemize}
We also neglect other constructions of algorithms.\\

\textbf{Algorithm 5.1}

Using the Euler's method~(\ref{eq31}) we solve numerically the last equation (\ref{eq09}) taking into account the initial
condition $y\left(x_0\right)=y_0$. Follow the discrete form of the Riemann-Liouville fractional integral~(\ref{eq30}) we
have
\begin{itemize}
 \item[\textbf{step 1}] prediction of necessary data: initial conditions, the fractional
               order $\alpha \in \langle 0, 1)$, the total length of
               calculations $x\in \langle x_0, x_N \rangle$, the step of
               calculations $h$.
 \item[\textbf{step 2}] Governing calculations: let $k = 1,\ldots, N$ then
               \begin{equation}
                \begin{array}{c} \label{eq41}
                 B_k  = h^{\alpha -1 } \bigg\{  - \lambda \Gamma \left( {2
                 -\alpha } \right)y_{k - 1} \\ 
                  - \sum\limits_{j = 2}^k {B_{j - 1}
                 \left[ {\left( {x_k  - x_{j - 2} } \right)^{1 - \alpha } –
                 \left( {x_k  - x_{j - 1} } \right)^{1 - \alpha } } \right]} 
                 \bigg\} \\~~\\ 
                 y_k = y_{k - 1}+ hB_k
                \end{array}
               \end{equation}

\end{itemize}

In summary we propose numerical schemes suitable for all possible cases of Eqn.~(\ref{eq42}). We then generally apply two
approaches. One is connected with the explicit method and the second is based on the predictor-corrector scheme.
It should be noted that we illustrate all the considered methods (the Euler's method, the Adams method and the Gear's method)
which are used in Eqns~(\ref{eq05}) (\ref{eq06}). However we limit this illustration in Eqn. (\ref{eq07})-(\ref{eq09})
only to the Euler's method. In this case we omit the predictor-corrector methods because there are problems with
applications. Using this approach we found that the Adams method requires known values of the first derivative in the
algebraic form. However the Gear's method requires the algebraic form of the second derivative which is used in
calculations of the correction factor.

\section{Results}

The algorithms presented in the previous section allow us to compare numerical and analytical results respectively. First we compare the results obtained from four discrete forms of the Caputo derivative. Thus we assume a~function
\begin{equation} \label{eq52}
 y\left( x \right) = x^2
\end{equation}
where the Caputo derivative is
\begin{equation} \label{eq53}
 g(x)={}_{x_0 }^C D_x^\alpha  x^2  = \frac{{\Gamma \left( 3 \right)x^{2 - \alpha
 } }}{{\Gamma \left( {3 - \alpha } \right)}}
\end{equation}
Table~\ref{tab01} \normalmarginpar \marginpar{\emph{Table~\ref{tab01}}} shows analytical values of the function~(\ref{eq53}) at assumed points. The other rows of this table present
the difference between the numerical and analytical results.
Analysing this table we can see that the liner-discrete form of the Caputo derivative (\ref{eq29}) (case-IV) gives
the best results. However, this discrete form is very complex. This is a~disadvantage of case-IV. It should be noted
that the middle-side form (\ref{eq28}) (case-III) also gives quite small errors. The next two discrete forms (\ref{eq26})
and (\ref{eq27}) generate large errors in comparison to the previous cases.

We try to estimate how several discrete forms of the Caputo derivative work in numerical schemes used to solve ordinary
differential equations. It should be remembered that, as presented in subsection~2.4, there are many constructions of
numerical approaches which are dependent on the form of the ordinary differential equation, the discrete form of the Caputo
derivative assumed (cases I to IV) and, finally, on the numerical scheme (the explicit or predictor-correct schemes)
used. Therefore, a~general question arises: how to efficiently compare all the possible constructions of numerical schemes
in order to comment on their practical application? First of all, we started our comparison between predictor-corrector
methods. We can apply three forms of the Caputo derivative (cases II to IV) and two numerical schemes (the Adams and
Gear's schemes) as presented in the previous section. Let us start with Eqn.~(\ref{eq05}) where initial conditions
\begin{equation} \label{eq55}
 y \left( 0 \right) = 0, ~~D^{1}y \left( 0 \right) = 1 
\end{equation}
are assumed. Note that Eqn.~(\ref{eq05}) has an analytical solution (\ref{eq17}).

\begin{figure}[!h]
\resizebox*{0.80\textwidth}{!}{\includegraphics{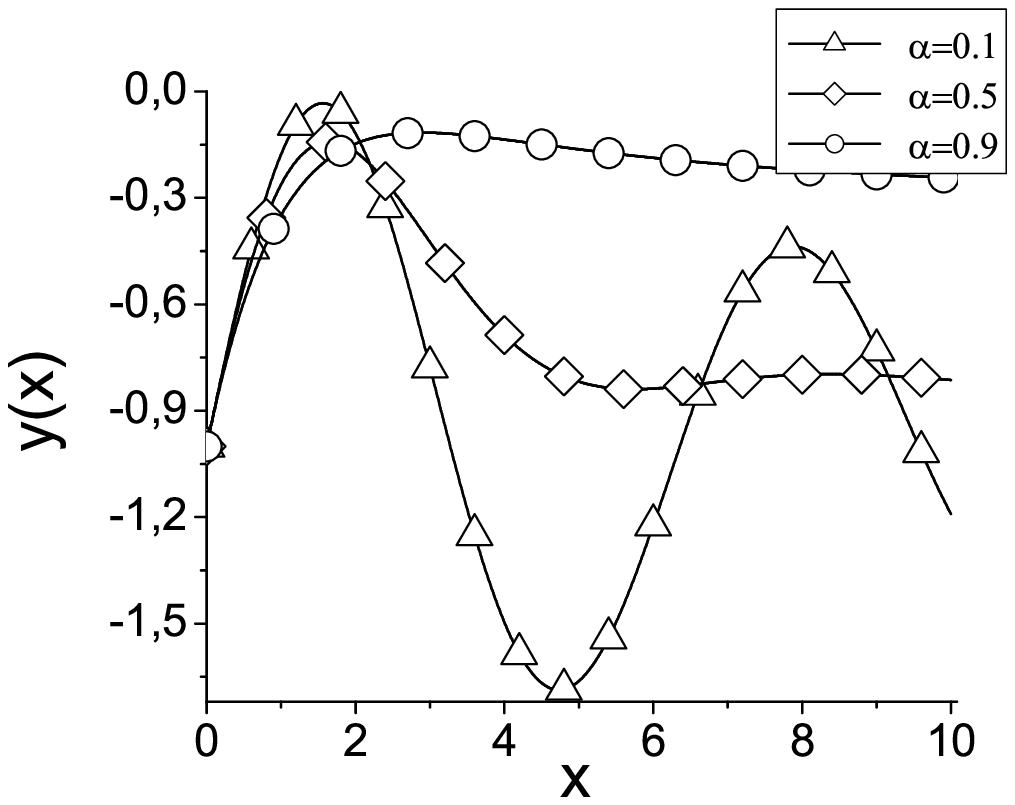}}\\
\resizebox*{0.80\textwidth}{!}{\includegraphics{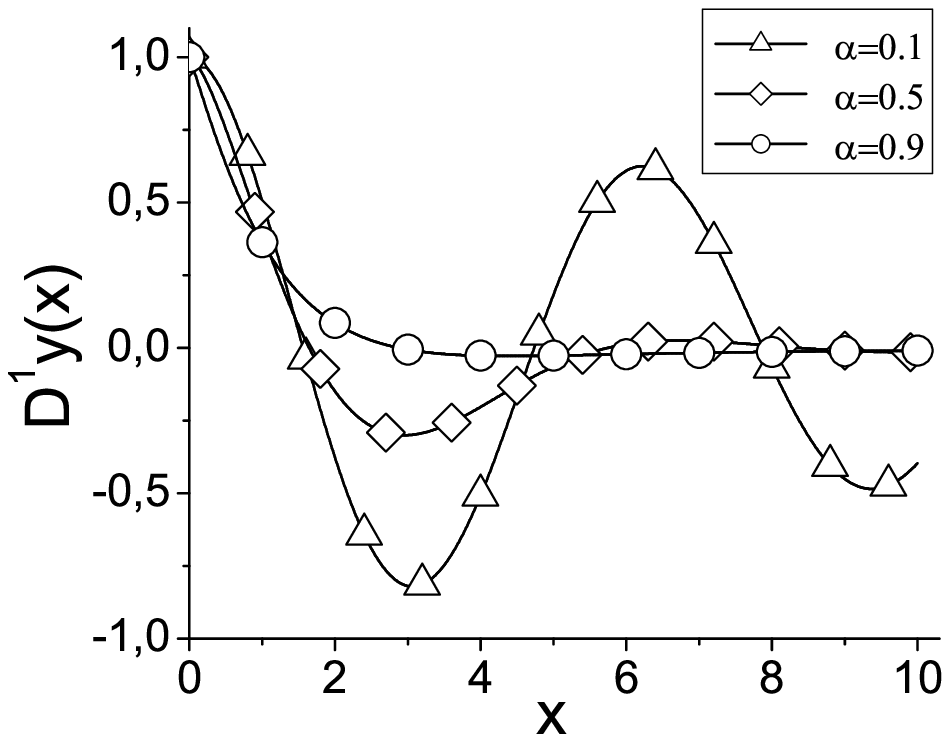}}
\caption{}\label{fig02} 
\end{figure}

Table~\ref{tab02} \reversemarginpar \marginpar{\emph{Table~\ref{tab02}}}
presents a direct comparison of the two predictor-corrector methods where three discrete forms of the Caputo derivative for each
method are included. 
It should be noted that we obtained many sets of data. Brief analysis presents the best results in the linear-discrete
form (case-IV) (\ref{eq29}) of the Caputo derivative for both the Gear's (\ref{eq34}) (\ref{eq35}) and Adams
(\ref{eq32}) (\ref{eq33}) methods. Satisfactory results can be also achieved by the application of the middle-side form
(case-III) (\ref{eq28}) obviously used for both methods. Taking into account the complexity of the forms, we chose the
middle-side form (\ref{eq28}) for the next calculations using the predictor-corrector method. We also noticed that the Gear's
method (\ref{eq34}) (\ref{eq35}) generates smaller error values than the Adams method (\ref{eq32}) (\ref{eq33}).
\begin{equation} \label{eq54}
 y \left( 0 \right) = -1, ~~D^{1}y \left( 0 \right) = 1 
\end{equation}
Now we consider all the numerical methods used to solve the same equation   (\ref{eq05}) with initial
conditions (\ref{eq54}). Fig.~\ref{fig02} \marginpar{\emph{Figure~\ref{fig02}}} illustrates the analytical solution taking into account the influence of  parameter $\alpha$. Notice that for such a solution we used algorithms 1.1, 1.2 and 1.3 respectively. Table~\ref{tab03} \marginpar{\emph{Table~\ref{tab03}}}
presents the analytical values given by formula (\ref{eq17}) and errors given by several algorithms used.
It can be observed that the Gear's method (algorithm 1.3) generates the smallest errors in comparison to other methods
used. However, the Euler's method (algorithm 1.1) performs better than the Adams method (algorithm 1.2). As we noted
in the previous section, predictor-corrector methods have some disadvantages which are revealed in the calculations of higher derivatives
for the Gear's method or necessary values of the function at the initial three points for the Adams method. Moreover,
a discrete form of the Caputo derivative has less of an influence on error calculations than a numerical method used in solving
an ordinary differential equation.

\begin{figure}[!h]
\resizebox*{0.80\textwidth}{!}{\includegraphics{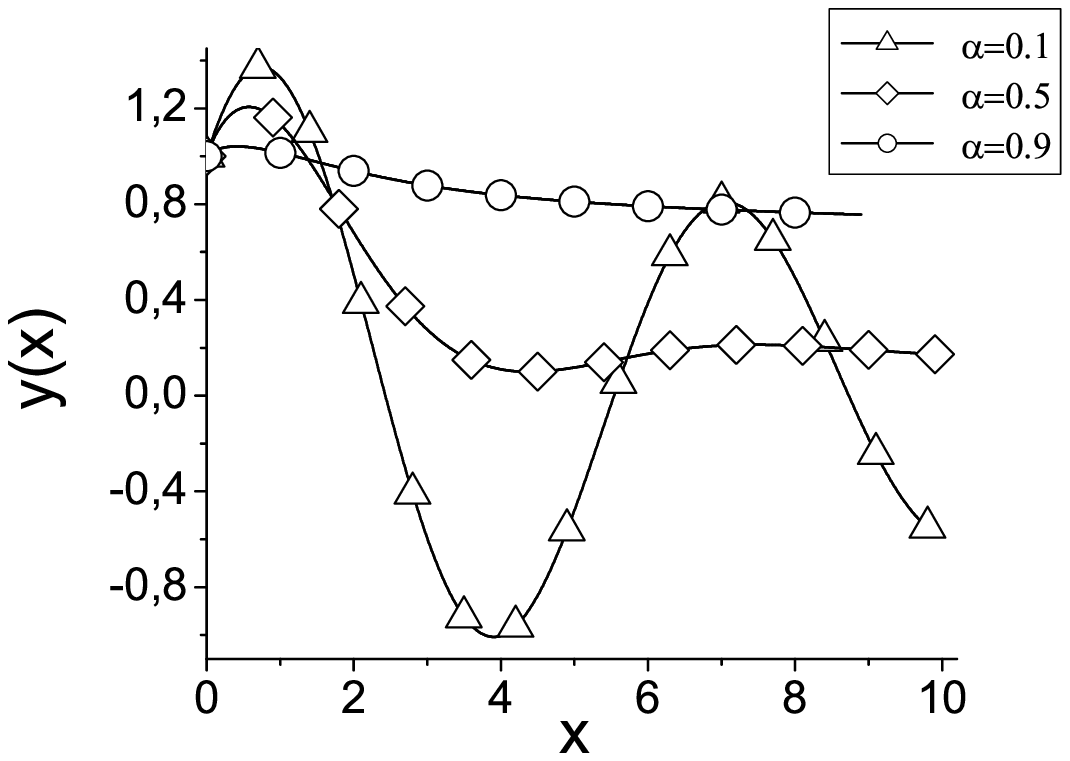}}
\resizebox*{0.80\textwidth}{!}{\includegraphics{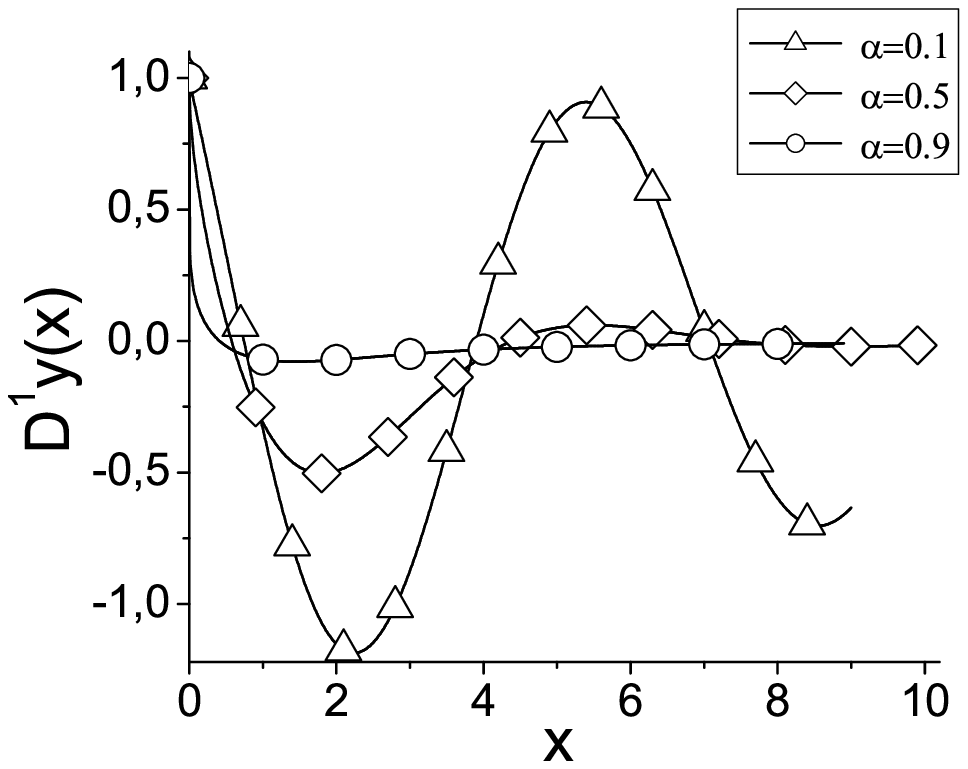}}
\caption{}\label{fig03}
\end{figure}

Next we consider all the numerical methods used to solve of Eqn.~(\ref{eq06}) with initial
conditions
\begin{equation} \label{eq57}
 y \left(0 \right) = 1,~~
 D^{1}y \left( 0 \right) = 1 
\end{equation}
This equation has an analytical solution given by formula (\ref{eq20})
It should be noted that for a such solution we used algorithms 2.1, 2.2 and 2.3 respectively. Fig.~\ref{fig03} \normalmarginpar \marginpar{\emph{Figure~\ref{fig03}}} shows
the behavior of the analytical solution over an independent value $x$ for different values of the parameter $\alpha$.
Table~\ref{tab04} \marginpar{\emph{Table~\ref{tab04}}} presents analytical values given by formula (\ref{eq20}) and errors given by several algorithms used.
This numerical solution confirms our previous considerations that the Euler's and Gear's methods show the smallest errors.
Taking into account both the errors generated by the numerical method and the complexity of numerical schemes we chose the Euler's
method as the method to use in the next calculations.

Assuming the initial condition $y\left( 0 \right)=2$ for Eqn. (\ref{eq07}) we obtain an analytical solution in the form
\begin{equation} \label{eq60}
 y \left( x \right)=2
\end{equation}

\begin{figure}[!h]
\resizebox*{0.90\textwidth}{!}{\includegraphics{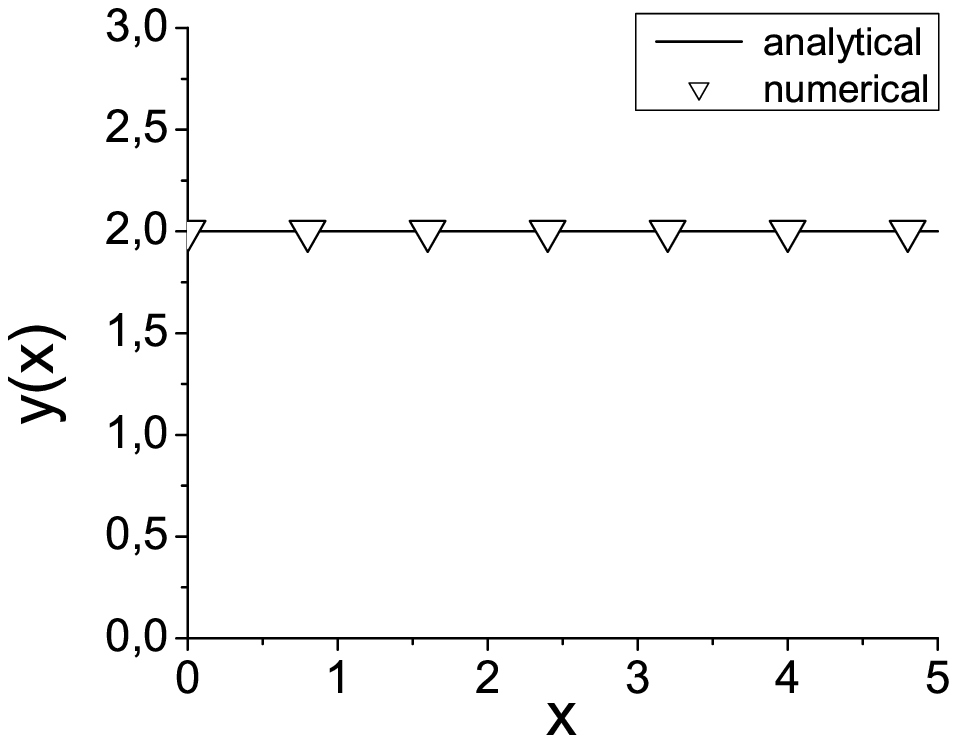}}
\caption{}\label{fig04}
\end{figure}

Fig.~\ref{fig04}  \marginpar{\emph{Figure~\ref{fig04}}} shows direct comparison between analytical solution (\ref{eq60}) and the Euler's method (algorithm 3.1) 
for Eqn. (\ref{eq07}). In this case we have not printed a~table because only a very small difference between the analytical and numerical
solutions can be observed.

\begin{figure}[!h]
\resizebox*{0.90\textwidth}{!}{\includegraphics{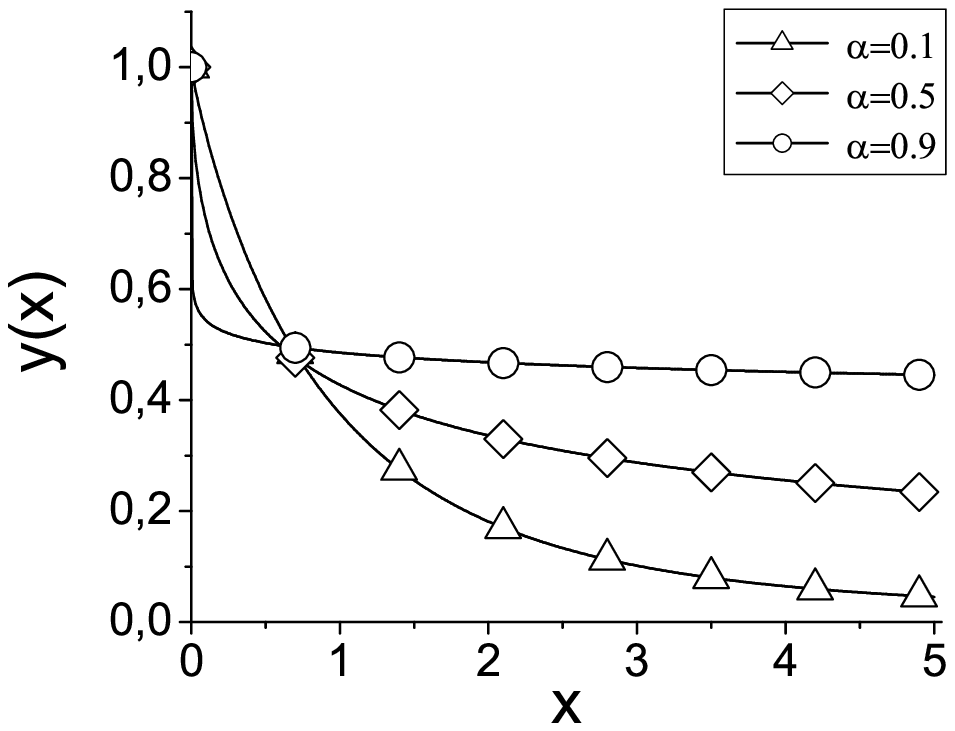}}
\caption{}\label{fig05}
\end{figure}

Next we consider Eqn.~(\ref{eq08}) with the initial condition $y(0)=1$. This equation has an analytical solution given by
formula~(\ref{eq23}) which is presented in Fig.~\ref{fig05}\marginpar{\emph{Figure~\ref{fig05}}}. We then apply three numerical schemes in order to compare data
with the analytical results. Table~\ref{tab05} \marginpar{\emph{Table~\ref{tab05}}} shows the analytical values given by formula (\ref{eq23}) and errors given by
the Euler's method (algorithm 4.1).
It may be observed that errors generated by the Euler's method are small. This confirms that the Euler's method is
an adequate method for solving ordinary differential equations with the a mixture of derivatives.

\begin{figure}[!h]
\resizebox*{0.90\textwidth}{!}{\includegraphics{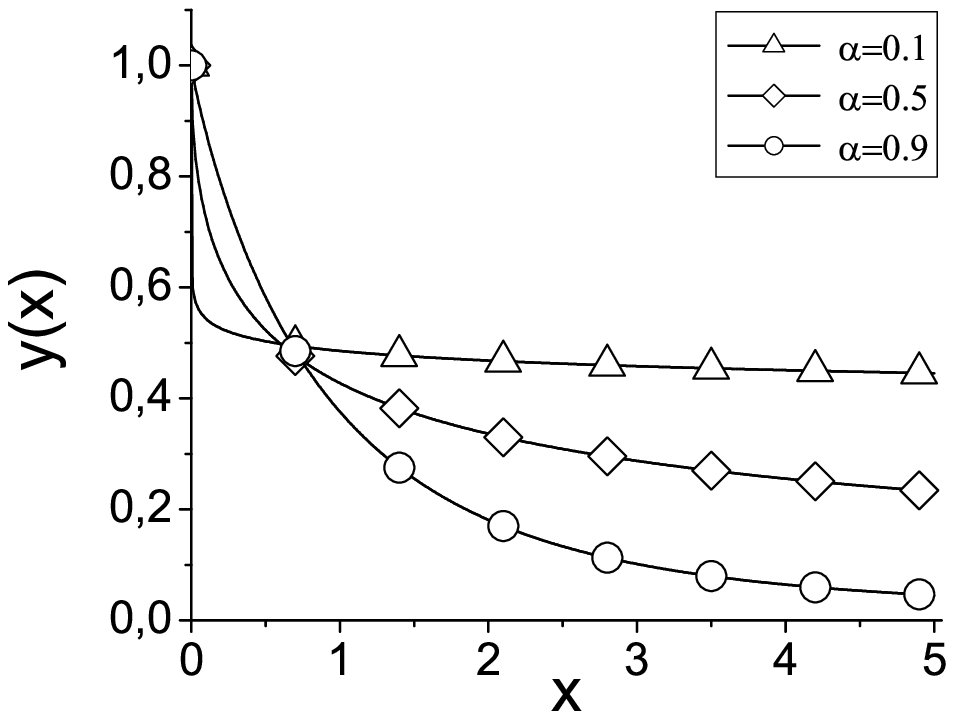}}
\caption{}\label{fig06}
\end{figure}

The last case concerns the fractional differential equation given by formula (\ref{eq09}). We assume the initial condition
$y(0)=1$. This equation has an analytical solution presented by formula (\ref{eq24}). We use the Euler's method given by algorithm 5.1 for a numerical solution. Fig.~\ref{fig06} \marginpar{\emph{Figure~\ref{fig06}}} presents the analytical solution (\ref{eq24}) over an
independent value $x$ for different values of parameter $\alpha$. Table~\ref{tab06} \marginpar{\emph{Table~\ref{tab06}}} shows analytical values and errors generated by the Euler's method (algorithm 5.1).
It can be observed that the Euler's method (algorithm 5.1) generates small error values. Summarising our results we proved
that the Euler's method is suitable for the numerical solution of ordinary differential equations having a mixture of
derivatives. This method has the following advantages: simplicity, a small difference between analytical and numerical
values depending on the step of calculations, stable error values over all the considered length of calculations.

\section{Conclusions}
In this study, we proposed numerical algorithms to solve ordinary differential equations where the
a mixture of fractional- and integer-order derivatives occurs. We used three known numerical techniques,
the Euler's, Adams and Gear's methods, to solve such equations. Taking into account the equation order,
we divided differential equations into three classes; where the integer order dominated over an integer
number calculated from the fractional order, where the integer order and number were the same and where the integer
number dominated over the integer order. In the considered equations we distinguished
two types of fractional operators: the left-side Riemann-Liouville and left-side Caputo derivatives. Using a known
transition rule between the derivatives, in a differential equation where the Riemman-Liouville operator occurs,
we changed the Riemann-Liouville derivative to the Caputo derivative. Next we proposed four discrete forms of the Caputo
derivative. On the basis of the previous classification of ordinary differential equations we illustrated the proper
algorithms. It should be noted that our algorithms are valid in the class of continuous functions. This assumption
allows us to solve the problem of how to include classical initial conditions into ordinary differential equations where
the Riemann-Liouville fractional derivative occurs.

Using direct comparison between the analytical and numerical data we obtained satisfactory results. Deeper analysis shows
that the \linebreak predictor-corrector methods (the Adams and Gear's methods) require the linear-discrete form (\ref{eq29})
of the Caputo derivative in order to reflect the analytical data more precisely. Satisfactory results can also be achieved by
the application of the middle-side form (\ref{eq28}) obviously used for both the methods. Comparing the results obtained
by the above methods we observed that the Gear's method gives better results than the Adams method. We also compare
the results obtained by the analytical solution, 
the predictor-corrector methods and the Euler's method which we first
successfully applied to solve ordinary differential equations including mixture of derivatives. We can say that
the Euler's and Gear's methods show the smallest errors. This may be unexpected because the Euler's method is a~method
of the order $O\left(h\right)$ and additionally it includes the left-side discrete form (\ref{eq26}) of the Caputo
derivative. On the other hand the predictor-corrector methods require additional circumstances, for example
the Adams method of fourth order needs the four initial values of the function 
to be determined and the Gear's method
needs higher derivatives in algebraic form at the starting point $x_0$. There are disadvantages of predictor-corrector
methods. Against this background we observed that a~discrete form of the Caputo derivative and the method order have less of an
influence on error calculations than the disadvantages of a~numerical method in solving an ordinary differential equation.
Taking into account both the errors generated by the numerical method and the complexity of numerical schemes we chose the Euler's
method as a~suitable method to use for practical calculations.

\newpage
List of captions for illustrations
\vspace{1cm}
\begin{itemize}
 \item[Fig. 1] Discrete forms of an integer derivative $D^n y(x)=B$ for the
               range $\langle x_{k-1},x_{k}\rangle$ 
               for $k=1,\ldots,N$:\protect\\
                 a) left-side\protect\\
                 b) right-side\protect\\
                 c) middle-side\protect\\
                 d) linear

\vspace{0.5cm}

 \item[Fig. 2]  Analytical solution of Eqn. (\ref{eq05}) with initial conditions
                (\ref{eq54}) over an independent value
                $x$ for different values of $\alpha$.

\vspace{0.5cm}

 \item[Fig. 3]  Analytical solution of Eqn. (\ref{eq06}) with initial conditions
                (\ref{eq57}) over an independent value
                $x$ for different values of $\alpha$.

\vspace{0.5cm}

 \item[Fig. 4]  Comparison of analytical (\ref{eq60}) and numerical
                (algorithm 3.1) for Eqn. (\ref{eq07}) with the
                initial condition $y(0)=2$.

\vspace{0.5cm}

 \item[Fig. 5]  Analytical solution of Eqn. (\ref{eq08}) with the initial
                condition \linebreak $y(0)=1$ over an independent value
                $x$ for different values of $\alpha$.

\vspace{0.5cm}

 \item[Fig. 6]  Analytical solution of Eqn. (\ref{eq09}) with the initial
                condition \linebreak $y(0)=1$ over an independent value
                $x$ for different values of $\alpha$.
\end{itemize}

\newpage

\newpage

\newpage

\newpage

\newpage

\newpage

\newpage
List of captions for tables
\vspace{1cm}

\begin{itemize}
 \item[Table~\ref{tab01}]  Analytical values of the Caputo
                           derivative (\ref{eq53}) and errors
                           generated by discrete forms of 
                           the Caputo derivative
                           (\ref{eq26}) (\ref{eq27}) (\ref{eq28}) 
                           (\ref{eq29}) for $h=0.01$.
                          
\vspace{0.5cm}

 \item[Table~\ref{tab02}]  Analytical solution of Eqn. (\ref{eq05}) with initial
                           conditions (\ref{eq55}) and errors generated by
                           the Gear's and Adams methods for $h=0.01$.

\vspace{0.5cm}

 \item[Table~\ref{tab03}]  Analytical results of Eqn.~(\ref{eq05}) with initial
                           conditions (\ref{eq54}) and errors generated by the
                           Euler's (algorithm~1.1), Adams
                           (algorithm~1.2) and Gear's (algorithm~1.3)
                           methods for $h=0.01$.

\vspace{0.5cm}

 \item[Table~\ref{tab04}]  Analytical results of Eqn.~(\ref{eq06}) with initial
                           conditions (\ref{eq57}) and errors generated by the
                           Euler's (algorithm~2.1), Adams
                           (algorithm~2.2) and Gear's (algorithm~2.3)
                           methods for $h=0.01$.

\vspace{0.5cm}

 \item[Table~\ref{tab05}]  Analytical results of Eqn.~(\ref{eq08}) with the
                           initial condition \linebreak $y(0)=1$ and errors
                           generated by the Euler's method (algorithm~4.1)
                           for $h=0.01$.

\vspace{0.5cm}

 \item[Table~\ref{tab06}]  Analytical results of Eqn.~(\ref{eq09}) with the
                           initial condition \linebreak $y(0)=1$ and errors 
                           generated by the Euler's method (algorithm~5.1)
                           for $h=0.01$.
\end{itemize}

\begin{table}
\caption{}\label{tab01}
\begin{tabular}{llllll}
\hline
{} & $y\left(1\right)$ & $y\left(4\right)$ & $y\left(6\right)$ & $y\left(8\right)$ & $y\left(10\right)$ \\
\hline
\hline
\multicolumn{6}{l}{$\alpha=0.1$}\\ \hline
analytical & 1.0944780 & 15.2447755 & 32.9377877 & 56.8955141 &  86.9374806\\
\hline
\multicolumn{6}{l}{errors generated by discrete forms of the Caputo
                     derivative}\\ \hline
case-I & 1.04e-2 & 3.62e-2 & 5.22e-2 & 6.76e-2 & 8.26e-2 \\
\hline
case-II & 1.04e-2 & 3.62e-2 & 5.21e-2 & 6.75e-2 & 8.26e-2 \\
\hline
case-III & 1.77e-5 & 1.98e-5 & 2.02e-5 & 2.07e-5 & 2.10e-5 \\\hline
case-IV & 0 & 0 & 2.00e-8 & 0 & 4.00e-8 \\
\hline\hline
\multicolumn{6}{l}{$\alpha=0.5$}\\ \hline
analytical & 1.5045056 & 12.0360444 & 22.1116256 & 34.0430746 & 47.5766431\\
\hline
\multicolumn{6}{l}{errors generated by discrete forms of the Caputo
                     derivative}\\ \hline
case-I & 1.17e-2 & 2.30e-2 & 2.81e-2 & 3.24e-2 & 3.61e-2 \\
\hline
case-II & 1.08e-2 & 2.21e-2 & 2.72e-2 & 3.14e-2 & 3.52e-2 \\
\hline
case-III & 4.60e-4 & 4.64e-4 & 4.65e-4 & 4.66e-4 & 4.66e-4 \\
\hline
case-IV & 0 & 0 & 0 & 2.00e-8 & 1.00e-8 \\
\hline
\hline
\multicolumn{6}{l}{$\alpha=0.9$}\\ \hline
analytical & 1.9111582 & 8.7813771 & 13.7171223 & 18.8232939 &  24.0600562\\
\hline
\multicolumn{6}{l}{errors generated by discrete forms of the Caputo
                     derivative}\\ \hline
case-I & 1.60e-2 & 1.76e-2 & 1.81e-2 & 1.85e-2 & 1.88e-2 \\
\hline
case-II & 4.98e-3 & 6.54e-3 & 7.04e-3 & 7.41e-3 & 7.70e-3 \\
\hline
case-III & 5.53e-3 & 5.53e-3 & 5.53e-3 & 5.53e-3 & 5.53e-3 \\
\hline
case-IV & 0 & 0 & 0 &1.00e-8& 1.00e-8 \\
\hline
\end{tabular}
\end{table}

\begin{table}
\caption{}\label{tab02}
\begin{tabular}{llllll}
\hline
{} & $y\left(1\right)$ & $y\left(4\right)$ & $y\left(6\right)$ & $y\left(8\right)$ & $y\left(10\right)$ \\
\hline
\hline
\multicolumn{6}{l}{$\alpha=0.1$}\\
\hline
analytical & 0.8226218 & -0.4960927 & -0.2222706 & 0.5587627 & -0.1911696\\
\hline
\multicolumn{6}{l}{errors generated by the Gear's method}\\
\hline
case-II & 2.71e-4 & 1.39e-2 & 4.21e-3 & 5.77e-3 & 1.39e-2 \\
\hline
case-III & 1.40e-5 & 1.24e-5 & 4.23e-6 & 1.47e-5 & 7.66e-6 \\
\hline
case-IV & 1.06e-5 & 1.75e-5 & 3.08e-5 & 5.57e-6 & 3.18e-5 \\
\hline
\multicolumn{6}{l}{errors generated by the Adams method}\\
\hline
case-II & 8.55e-5 & 1.40e-2 & 4.27e-3 & 5.89e-3 & 1.40e-2 \\
\hline
case-III & 1.71e-4 & 1.05e-4 & 5.35e-5 & 1.11e-4 & 3.97e-5 \\
\hline
case-IV & 1.70e-4 & 1.14e-4 & 3.67e-5 & 1.09e-4 & 5.85e-5 \\
\hline
\hline
\multicolumn{6}{l}{$\alpha=0.5$}\\ 
\hline
analytical & 0.7374822 & 0.3129516 & 0.1623955 & 0.2017797 &  0.1867275\\
\hline
\multicolumn{6}{l}{errors generated by the Gear's method}\\
\hline
case-II & 7.85e-4 & 7.84e-3 & 5.76e-3 & 4.45e-3 & 5.04e-3 \\
\hline
case-III & 1.26e-4 & 1.02e-4 & 1.96e-5 & 2.53e-5 & 8.40e-6 \\
\hline
case-IV & 1.67e-4 & 7.25e-5 & 4.27e-5 & 4.71e-5 & 4.29e-5 \\
\hline
\multicolumn{6}{l}{error generated by the Adams method}\\
\hline
case-II & 6.92e-5 & 7.54e-3 & 5.62e-3 & 4.27e-3 & 4.87e-3 \\
\hline
case-III & 5.90e-4 & 3.97e-4 & 1.22e-4 & 1.56e-4 & 1.75e-4 \\
\hline
case-IV & 5.40e-4 & 2.22e-4 & 1.06e-4 & 1.34e-4 & 1.23e-4 \\
\hline 
\hline
\multicolumn{6}{l}{$\alpha=0.9$}\\
\hline
analytical & 0.6512921 & 0.8643086 & 0.8128471 & 0.7787088 &  0.7563200\\
\hline
\multicolumn{6}{l}{errors generated by the Gear's method}\\
\hline
case-II & 2.29e-3 & 4.97e-3 & 5.13e-3 & 5.09e-3 & 5.06e-3 \\
\hline
case-III & 1.09e-3 & 5.18e-5 & 1.20e-4 & 7.50e-5 & 4.63e-5 \\
\hline
case-IV & 1.80e-3 & 2.39e-3 & 2.25e-3 & 2.15e-3 & 2.09e-3 \\
\hline
\multicolumn{6}{l}{errors generated by the Adams method}\\
\hline
case-II & 8.76e-5 & 1.78e-3 & 2.13e-3 & 2.23e-3 & 2.28e-3 \\
\hline
case-III & 1.25e-3 & 3.19e-3 & 3.06e-3 & 2.89e-3 & 2.78e-3 \\
\hline
case-IV & 5.65e-4 & 7.89e-4 & 7.37e-4 & 7.02e-4 & 6.79e-4 \\
\hline
\end{tabular}
\end{table}

\begin{table}
\caption{}\label{tab03}
\begin{tabular}{llllll}
\hline
{} & $y\left(1\right)$ & $y\left(4\right)$ & $y\left(6\right)$ & $y\left(8\right)$ & $y\left(10\right)$ \\
\hline
\hline
\multicolumn{6}{l}{$\alpha=0.1$}\\ \hline
analytical & -0.1773782 & -1.4960927 & -1.2222706 & -0.4412373 & -1.1911696\\
\hline
Euler & 1.79e-5 & 2.68e-5 & 2.92e-5 & 4.23e-6 & 3.85e-5 \\
\hline
Gear & 1.40e-5 & 1.24e-5 & 4.18e-6 & 1.47e-5 & 7.60e-6\\
\hline
Adams & 1.70e-4 & 1.14e-4 & 3.67e-5 & 1.09e-4 & 5.85e-5 \\
\hline\hline
\multicolumn{6}{l}{$\alpha=0.5$}\\ \hline 
analytical & -0.2625177 & -0.6870484 & -0.8376044 & -0.7982203 & -0.8132725\\
\hline
Euler & 1.26e-4 & 1.04e-4 & 2.41e-5 & 2.51e-5 & 9.32e-6 \\
\hline
Gear & 1.26e-4 & 1.02e-4 & 1.96e-5 & 2.53e-5 & 8.40e-6 \\
\hline
Adams &  5.40e-4 & 2.22e-4 & 1.06e-4 & 1.34e-4 & 1.23e-4 \\
\hline\hline
\multicolumn{6}{l}{$\alpha=0.9$}\\ \hline
analytical & -0.3487079 & -0.1356914 & -0.1871529 & -0.2212912 &  -0.2436800\\
\hline
Euler &1.08e-3 & 5.29e-5 & 1.20e-4 & 7.49e-5 & 4.62e-5 \\
\hline
Gear & 1.06e-3 & 5.18e-5 & 1.20e-4 & 7.50e-5 & 4.63e-5\\
\hline
Adams & 5.65e-4 & 7.89e-4 & 7.37e-4 & 7.02e-4 & 6.79e-4 \\ \hline
\end{tabular}
\end{table}

\begin{table}
\caption{}\label{tab04}
\begin{tabular}{llllll}
\hline
{} & $y\left(1\right)$ & $y\left(4\right)$ & $y\left(6\right)$ & $y\left(8\right)$ & $y\left(10\right)$ \\
\hline
\hline
\multicolumn{6}{l}{$\alpha=0.1$}\\ \hline
analytical & 1.3290813 & -1.0029826 & 0.3872104& 0.4928008 & -0.5873420\\
\hline
Euler & 4.83e-3 & 3.40e-3 & 8.51e-4 & 3.06e-3 & 1.49e-3 \\
\hline
Gear & 4.83e-3 & 3.41e-3 & 8.92e-4 & 3.10e-3 & 1.49e-3\\
\hline
Adams & 1.06e-2 & 6.39e-3 & 2.50e-3 & 7.28e-3 & 2.39e-3 \\
\hline
\hline
\multicolumn{6}{l}{$\alpha=0.5$}\\ \hline
analytical & 1.1341116 & 0.1100800 & 0.1748923 & 0.2090891 & 0.1714270\\
\hline
Euler & 7.42e-3 & 8.06e-5 & 3.56e-4 & 8.41e-4 & 6.45e-4 \\
\hline
Gear & 7.42e-3 & 7.54e-5 & 3.53e-4 & 8.43e-4 & 6.45e-4\\
\hline
Adams & 2.27e-2 & 7.88e-3 & 4.71e-3 & 6.06e-3 & 5.61e-3 \\
\hline\hline
\multicolumn{6}{l}{$\alpha=0.9$}\\ \hline
analytical & 1.0146791 & 0.8374876 & 0.7913672 & 0.7652525 & -\\
\hline
Euler & 5.73e-3 & 2.17e-3 & 2.01e-3 & 1.96e-3 & - \\
\hline
Gear & 5.74e-3 & 2.17e-3 & 2.01e-3 & 1.96e-3 & -\\
\hline
Adams & 3.15e-2 & 3.69e-2 & 3.50e-2 & 3.37e-2 & - \\
\hline
\end{tabular}
\end{table}

\begin{table}
\caption{}\label{tab05}
\begin{tabular}{llllll}
\hline
{} & $y\left(1\right)$ & $y\left(2\right)$ & $y\left(3\right)$ & $y\left(4\right)$ & $y\left(5\right)$ \\
\hline
\hline
\multicolumn{6}{l}{$\alpha=0.1$}\\
\hline
analytical & 0.3760660 & 0.1811155 & 0.1014866 & 0.0644356 & 0.0452231\\
\hline
Euler & 1.36e-3 & 9.40e-4 & 5.42e-4 & 3.03e-4 & 1.70e-4 \\
\hline
\hline
\multicolumn{6}{l}{$\alpha=0.5$}\\
\hline
analytical & 0.4275836 & 0.3362040 & 0.2873412 & 0.2553957 & 0.2323262\\
\hline
Euler & 9.61e-4 & 8.60e-4 & 7.88e-4 & 7.31e-4 & 6.85e-4 \\
\hline
\hline
\multicolumn{6}{l}{$\alpha=0.9$}\\
\hline
analytical & 0.4855645 & 0.4682030 & 0.4580801 & 0.4509182 & 0.4453768\\
\hline
Euler & 1.23e-3 & 1.19e-3 & 1.17e-3 & 1.16e-3 & 1.14e-3 \\
\hline
\end{tabular}
\end{table}

\begin{table}
\caption{}\label{tab06}
\begin{tabular}{llllll}
\hline
{} & $y\left(1\right)$ & $y\left(2\right)$ & $y\left(3\right)$ & $y\left(4\right)$ & $y\left(5\right)$ \\
\hline
\hline
\multicolumn{6}{l}{$\alpha=0.1$}\\
\hline
analytical & 0.4855645 & 0.4682030 &  0.4580801 & 0.4509182 & 0.4453768\\
\hline
Euler & 1.33e-4 & 7.06e-5 & 4.85e-5 & 3.69e-5 & 3.07e-5 \\
\hline
\hline
\multicolumn{6}{l}{$\alpha=0.5$}\\
\hline
analytical & 0.4275836 & 0.3362040 & 0.2873412 & 0.2553957 & 0.2323262\\
\hline
Euler & 8.39e-4 & 5.30e-4 & 3.73e-4 & 2.82e-4 & 2.23e-4 \\
\hline
\hline
\multicolumn{6}{l}{$\alpha=0.9$}\\
\hline
analytical & 0.3760660 & 0.1811155 & 0.1014866 & 0.0644356 & 0.0452231\\
\hline
Euler & 1.55e-3 & 1.27e-3 & 7.99e-4 & 4.77e-4 & 2.88e-4 \\
\hline
\end{tabular}
\end{table}

\end{document}